\input amstex\documentstyle{amsppt}  
\pagewidth{12.5cm}\pageheight{19cm}\magnification\magstep1  
\topmatter
\title Strata of a disconnected reductive group\endtitle
\author G. Lusztig\endauthor
\address{Department of Mathematics, M.I.T., Cambridge, MA 02139}\endaddress
\thanks{Supported by NSF grant DMS-1855773.}\endthanks
\endtopmatter   
\document
\define\hBP{\widehat{BP}}
\define\hg{\hat g}

\define\uR{\un R}

\define\rk{\text{\rm rk}}

\define\Irr{\text{\rm Irr}}

\define\si{\sim}

\define\sqc{\sqcup}

\define\op{\oplus}
   
\define\part{\partial}
\define\emp{\emptyset}

\define\n{\notin}

\define\m{\mapsto}
\define\do{\dots}

\define\lra{\leftrightarrow}

\define\sub{\subset}    

\define\T{\times}
\define\ti{\tilde}
\define\nl{\newline}
\redefine\i{^{-1}}

\define\un{\underline}

\define\ot{\otimes}
\define\bbq{\bar{\QQ}_l}

\define\Ad{\text{\rm Ad}}
\define\Hom{\text{\rm Hom}}

\define\sgn{\text{\rm sgn}}

\define\supp{\text{\rm supp}}

\define\a{\alpha}
\redefine\b{\beta}

\define\g{\gamma}
\redefine\d{\delta}
\define\e{\epsilon}

\define\io{\iota}
\redefine\o{\omega}
\define\p{\pi}

\redefine\t{\tau}

\define\k{\kappa}
\redefine\l{\lambda}
\define\z{\zeta}
\define\x{\xi}

\define\vt{\vartheta}

\redefine\G{\Gamma}

\define\Si{\Sigma}

\redefine\L{\Lambda}

\define\boc{\bold c}

\define\kk{\bold k}

\define\CC{\bold C}

\define\EE{\bold E}

\define\NN{\bold N}

\define\QQ{\bold Q}

\define\ZZ{\bold Z}

\define\ca{\Cal A}
\define\cb{\Cal B}
\define\cc{\Cal C}

\define\cg{\Cal G}

\define\ci{\Cal I}

\define\cp{\Cal P}

\define\cs{\Cal S}
\define\ct{\Cal T}

\define\cz{\Cal Z}
\define\cx{\Cal X}
\define\cy{\Cal Y}
\define\ucx{\un{\cx}}
\define\ucy{\un{\cy}}

\define\te{\ti e}

\define\tI{\ti I}

\define\bI{\bar I}

\define\che{\check}

\define\tcb{\ti{\cb}}

\head Introduction\endhead
\subhead 0.1\endsubhead
Let $\kk$ be an algebraically closed field of characteristic $p\ge0$ and let $G$ be a possibly disconnected
reductive algebraic group over $\kk$ with a fixed connected component $D$. The identity component $G^0$ of $G$
acts on $D$ by conjugation.
In the case where $D=G^0$, Steinberg \cite{St65} has defined the open set $D_{reg}$ of regular elements in
$G^0$ (a Lie algebra analogue of this set was earlier defined by Kostant); in \cite{L15} we have defined
a partition of $G^0$ into finitely many Strata, one of which is $D_{reg}$. The goal of this paper is to define
(without assuming $D=G^0$) a partition of $D$ into finitely many Strata, (each of which is a union of
$G^0$-conjugacy classes of fixed dimension) generalizing the partition \cite{L15}  of $G^0$ into Strata.
(We use the term ``Stratum'' of $D$ to distinguish it from ``stratum'' of $D$ in \cite{L03,\S3}. In fact,
every Stratum of $D$ is a finite union of strata of $D$.)
To define the Strata of $D$
we associate to any $g\in D$ an irreducible representation $E_g$ (a variant of Springer's
representation) of the subgroup of the Weyl group given by the fixed point set of the action
of $D$ on the Weyl
group. Then we say that $g,g'$ in $D$ are in the same Stratum whenever $E_g,E_{g'}$ are isomorphic.
We show that the Strata of $D$ are indexed by a set defined purely in terms of the Weyl group and
its automorphism defined by $D$ (thus extending a result of \cite{L15}).
A definition of Strata of $D$ (different from the one in this paper) was sketched without proof
(under some additional assumptions on $G,D$) in \cite{L15, 6.1}. 

I wish to thank the referee for very helpful comments.

\subhead 0.2\endsubhead
{\it Notation.} We fix a prime number $l$ invertible in $\kk$. Let $\bbq$ be an algebraic closure of the
field of $l$-adic numbers. For
any finite group $\G$ let $Irr(\G)$ be the set of irreducible representations (over $\bbq$) of $\G$ (up to
isomorphism). For a Weyl group $W$ and for $E\in\Irr(W)$ we denote by $b_E$ the smallest integer $\ge0$ such that
$E$ appears in the $b_E$-th symmetric power of the reflection representation of $W$.

\head Contents \endhead
1. Definition of Strata.

2. Bipartitions.

3. Examples.

4. Proof of Theorem 1.15.

5. Further results.

\head 1. Definition of Strata\endhead
\subhead 1.1\endsubhead
Let $\cb=\cb_{G^0}$ be the flag manifold of $G^0$. Let $\ct=\ct_{G^0}$ be the set of maximal tori of $G^0$.
Let $\tcb=\tcb_{G^0}$ be the set of all pairs $(T,B)$ where $B\in\cb,T\in\ct,T\sub B$. Let $W=W_{G^0}$ be the
Weyl group of $G^0$, viewed as the set of $G^0$-orbits on $\cb\T\cb$; this is naturally a (finite) Coxeter group
with length function $w\m|w|$. For $(B_1,B_2)\in\cb\T\cb$ let $pos(B_1,B_2)\in W$ be the $G^0$-orbit of
$(B_1,B_2)$. (The product of $pos(B_1,B_2)$ with $pos(B_2,B_3)$ is $pos(B_1,B_3)$ provided that some maximal
torus of $G^0$ is contained in $B_1,B_2,B_3$.)

Now $G$ acts by inner automorphisms on $G^0$ and this induces an action of $G/G^0$ on
$W$. In particular, $D$, viewed as an element of $G/G^0$,
defines an automorphism $[D]:W@>>>W$ whose fixed point set is denoted by $W^D$. For each $[D]$-orbit $o$ on
$\{s\in W;|s|=1\}$ we denote by $s_o$ the element of maximal length in the subgroup of $W$ generated by $o$;
then the elements $s_o\in W^D$ for various $o$ are the simple reflections for a Coxeter group structure on $W^D$. (See \cite{L14, Appendix}.)
For $g\in G$ let $g_s$ be the semisimple part of $g$ and let $g_u$ be the unipotent part of $g$.

\subhead 1.2\endsubhead
For $g\in G$ let $\tcb_g=\{(T,B)\in\tcb;gTg\i=T,gBg\i=B\}$. Let $D_{un}:=\{g\in D;g=g_u\}$. For any
$e=(T,B)\in\tcb$ let $S=S_{e;D}$ be the set of all $g\in D$ such that $e\in\tcb_g$.
This is a single orbit of $T$ acting on $D$ by left multiplication  
(resp. by right multiplication). Now $T$ is uniquely determined by $S$ (it is the set of all
$s's\i$ where $(s,s')\in S\T S$). Let ${}^*S_{e;D}=\{g\in S;Z_G(g_s)^0\sub T\}$; this is an open dense
subset of $S$, see \cite{L04, 3.11}. (Here $Z_{-}(-)$ denotes a centralizer).
Let ${}^*D=\cup_{e\in\tcb}{}^*S_{e;D}$. This is a stratum of $G$ in the
sense of \cite{L04, 3.1}. By \cite{L04, 3.13, 3.16}, ${}^*D$ is an irreducible locally closed subset of $D$
of dimension $\dim(G^0/T)+\dim S=\dim(G^0)=\dim D$. Hence ${}^*D$ is an open dense subset of $D$.
In the case where $D=G^0$, ${}^*D$ is the set of regular semisimple elements in $G^0$.
(In the general case, ${}^*D$ is the same as the subset of $D$ associated in \cite{L87, (2.3.2)} to $(L^0,\Si)$
with $(L^0,\Si)$ as in \cite{L87, (2.3.1)} with $L^0\in\ct$ and with the added
requirement \cite{L87, (2.3.5)} which, contrary to what is stated in {\it loc.cit.}, is not an automatic
consequence of the conditions \cite{L87, (2.3.1)}.)
Now ${}^*D$ has a finite unramified covering ${}^*\p:{}^*\ti D@>>>{}^*D$, see \cite{L04, 3.13} or
\cite{L87, 2.5}. In
our case we have ${}^*\ti D=\{(g,B')\in{}^*D\T\cb;gB'g\i=B'\}$, ${}^*\p(g,B')=g$. If $g\in {}^*S_{e;D}$
with $e=(T,B)\in\tcb$, then ${}^*\p\i(g)$ is the finite set
$\{(g,B');B'\in\cb,T\sub B',pos(B,B')\in W^D\}$. As stated in {\it loc.cit.}, ${}^*\p$ is a principal
covering 
whose group is in our case $W^D$, which acts freely on ${}^*\ti D$
by $w:(g,B')\m(g,B'')$ where $B''\in\cb$ is defined by $T\sub B''$, $pos(B',B'')=w$.
Let $\ti D=\{(g,B')\in D\T\cb;gB'g\i=B'\}$. Define $\p:\ti D@>>>D$ by $\p(g,B')=g$. This is a proper
surjective morphism. According to \cite{L87, 2.6}, \cite{L04} (see also  \cite{S04} in the case where
$D_{un}\ne\emp$), this morphism is
small. (When $D=G^0$, this is an observation of \cite{L81}.) It follows that $\p_!(\bbq)$ is
the intersection cohomology complex of $D$ with coefficients in the local system ${}^*\p_!(\bbq)$ on
${}^*D$, which has a natural action of $W^D$. 
Then $\p_!(\bbq)[\dim D]=\op_EE\ot\p_!(\bbq)[\dim D]_E$ where $E$ runs through $\Irr(W^D)$,
$\p_!(\bbq)[\dim D]_E$ is a simple perverse sheaf on $D$ and the action of $W^D$ on ${}^*\p_!(\bbq)$ restricts
for each $E$ to the obvious $W^D$ action on $E$ tensor the identity on $\p_!(\bbq)[\dim D]_E$.

\subhead 1.3\endsubhead
In this subsection we assume that $D_{un}\ne\emp$.
Let $g\in D_{un}$ and let $\boc\sub D_{un}$ be the $G^0$-conjugacy class of $g$.

(a) {\it There is a unique $E=E_g\in\Irr(W^D)$ such that $(\p_!(\bbq)[\dim D]_E)|D_{un}$ is (up to shift) the
intersection cohomology complex of the closure of $\boc$ with coefficients in $\bbq$.}
\nl
In the case where $D=G^0$ this is proved in \cite{L84a}. Similar arguments apply in the general case, see
\cite{L04, 8.2(b)}, \cite{S04}.

\subhead 1.4\endsubhead
For any $g\in G$ let $\cb_g=\{B\in\cb;gBg\i=B\}$; this is a closed nonempty subvariety of $\cb$, see
\cite{St68}. Let $\un\cb_g$ be the set of irreducible components of $\cb_g$.

Recall that an element $h\in G$ is said to be quasi-semisimple (qss) if $\tcb_h\ne\emp$. In this subsection we
fix a qss element $h$ of $G$; let $H=Z_G(h)$.

(a) {\it $H$ is reductive, see \cite{Sp82, II, 1.17}, and $(H\cap G^0)/H^0$ consists of semisimple elements, see
\cite{DM94, 1.8};}

(b) {\it if $B\in\cb_h$, then $B\cap H^0\in\cb_{H^0}$; for any $Z\in\un\cb_h$, $B\m B\cap H^0$ is an isomorphism
$t_Z:Z@>\si>>\cb_{H^0}$, see \cite{Sp82, II, 1.17}.}

(c) {\it if $T'\in\ct_{H^0}$ then $T:=Z_{G^0}(T')\in\ct$; moreover $T$ is the unique maximal torus of
$G^0$ that contains $T'$, see  \cite{Sp82, II, 1.15}, \cite{DM94, 1.8};}

(d) {\it if $T\in\ct$ is such that $hTh\i=T,hBh\i=B$ for some $B\in\cb$ containing $T$ then
$T':=T\cap H^0$ is in $\ct_{H^0}$, see \cite{Sp82, II, 1.15}, \cite{DM94, 1.8}.}

\subhead 1.5\endsubhead
Let $h,T',T,H$ be as in 1.4(c). If $n\in H^0,nT'n\i=T'$, then clearly $nTn\i=T$. Thus we have a well
defined inclusion $N_{H^0}(T')\sub N_{G^0}(T)$ (here $N_{-}(-)$ denotes a normalizer); this carries $T'$
into $T$ hence induces a homomorphism

(a) $N_{H^0}(T')/T'@>>> N_{G^0}(T)/T$.
\nl
This homomorphism is injective. It is enough to show that $N_{H^0}(T')\cap T\sub T'$; this follows from
$H^0\cap T=T'$ (see 1.4).
Now $\Ad(h)$ defines an automorphism of each side of (a) and the map (a) is compatible with these
automorphisms; moreover the automorphism on the left hand side of (a) is the identity. It follows that (a)
restricts to

(b) {\it an imbedding $N_{H^0}(T')/T'@>>> (N_{G^0}(T)/T)^h$}
\nl
where ${}^h$ denotes the fixed point set of $\Ad(h)$.

\subhead 1.6\endsubhead
Let $g\in D$. Note that $g_s$ is qss in $G$ (see \cite{St68}). Hence the results in 1.4, 1.5 apply with
$h=g_s$. Let $G'=Z_G(g_s)$. Let $D'$ be the connected component of $G'$ that
contains $g$. Let $\cb'=\cb_{G'{}^0}$, $W'=W_{G'{}^0}$. We define $[D']:W'@>>>W'$ and $W'{}^{D'}$ in terms of
$W',D'$ in the same way as $[D]:W@>>>W$ and $W^D$ were defined in terms of $W,D$.
For any $Z\in\un\cb_{g_s}$, the assignment

$G'{}^0-\text{orbit of }(B'_1,B'_2)\m pos(t_Z\i(B'_1),t_Z\i(B'_2))$
\nl
is a map $\t_Z:W'@>>>W$.

Let $B'\in\cb'$ and let $T'$ be a maximal torus of $B'$. Let $B=t_Z\i(B')\in Z$, and let $T$ be the
unique maximal torus of $G^0$ that contains $T'$; we have $T\sub B$. We have a commutative diagram
$$\CD  N_{G'{}^0}(T')/T'@>j>> N_{G^0}(T)/T \\
        @Va'VV                  @VaVV        \\
        W'@>\t_Z>>               W          \\
        \endCD$$
where $a'(nT')=G'{}^0-\text{orbit of }(B',nB'n\i)$, $a(nT)=G^0-\text{orbit of }(B,nBn\i)$ and $j$ is as in
1.5(a) with $h=g_s$, $H=G'$. (We use that for $n\in N_{G'{}^0}(T')$ we have $t_Z(nBn\i)=nB'n\i$.) Since $j$
is an injective homomorphism it follows that $\t_Z$ is an injective homomorphism.

Now let $Z,Z'$ in $\un\cb_{g_s}$; let $w=w_{Z,Z'}$ be the unique $G^0$-orbit on $\cb\T\cb$ that
contains $(t_Z\i(B'),t_{Z'}\i(B'))$ for any $B'\in\cb'$ (we view $w$ as an element of $W$). We show that

(a) $\t_{Z'}(w')w=w\t_Z(w')$ for any $w'\in W'$.
\nl
It is enough to show that if $B'_1,B'_2$ are in $\cb'$ then $t_Z\i(B'_1),t_Z\i(B'_2),t_{Z'}\i(B'_2)$ contain
a common maximal torus and that $t_Z\i(B'_1),t_{Z'}\i(B'_1),t_{Z'}\i(B'_2)$ contain a common maximal torus.
Now let $T'$ be a maximal torus of $G'{}^0$ contained in $B'_1\cap B'_2$.
Let $T$ be the unique maximal torus of $G^0$ such that $T'\sub T$, see 1.4(c).
We have $T\sub t_Z\i(B'_1),T\sub t_Z\i(B'_2),T\sub t_{Z'}\i(B'_1),T\sub t_{Z'}\i(B'_2)$. Our claim is proved:
the orbit represented by either side of (a) is that of $(t_Z\i(B'_1),t_{Z'}\i(B'_2))$.

Now $G'$ acts (by conjugation) on $\cb_{g_s}$; this induces an action $g':Z\m g'(Z)$
of $G'/G'{}^0$ on $\un\cb_{g_s}$. From the definition, for any $Z\in\un\cb_{g_s}$ and any $g'\in G'/G'{}^0$
we have $g'\t_Z=\t_{g'(Z)}g'$ as maps $W'@>>>W$. In particular, taking $g'=g$ we obtain
$[D](\t_Z(w'))=\t_{g(Z)}([D'](w'))$ for any $w'\in W'$. Let $(\un\cb_{g_s})_g=\{Z\in\un\cb_{g_s};g(Z)=Z\}$.
This set is nonempty; indeed, for some $B\in\cb$ we have $gBg\i=B$ (see \cite{St68, p.49}) so that we have also
$g_sBg_s\i=B$ and if $Z_0$ is the irreducible component of $\cb_{g_s}$ that contains $B$ then
$Z_0\in(\un\cb_{g_s})_g$. Note that for any $Z\in(\un\cb_{g_s})_g$ we have $[D]\t_Z=\t_Z[D']$ as maps
$W'@>>>W$; hence $\t_Z$ restricts to an (injective) homomorphism $\t_Z^D:W'{}^{D'}@>>>W^D$.

Assume now that  $Z,Z'$ are in $(\un\cb_{g_s})_g$; let $w=w_{Z,Z'}$ be as above. We show that

(b) $w\in W^D$.
\nl
It is enough to show that if $B'\in\cb'$ then

$(gt_Z\i(B')g\i,gt_{Z'}\i(B')g\i)=(t_Z\i(gB'g\i),t_{Z'}\i(gB'g\i))$
\nl
is in the same $G^0$-orbit as $(t_Z\i(B'),t_{Z'}\i(B'))$. This is clear since $gB'g\i,B'$ are Borel subgroups
of $G'{}^0$ hence are in the same $G'{}^0$-orbit.

From (a),(b) we see that

(c) {\it the injective homomorphism $\t_Z^D:W'{}^{D'}@>>>W^D$ defined for $Z\in(\un\cb_{g_s})_g$ is
independent of the choice of $Z$, up to composition with an inner automorphism of $W^D$.}

\subhead 1.7\endsubhead
Let $h\in D$. Assume that $h$ is qss in $G$. Then the results in 1.4, 1.5 apply to $h$ and $H=Z_G(h)$. Let
$\cb'=\cb_{H^0}$, $W'=W_{H^0}$. For any $Z\in\un\cb_h$, the assignment

$H^0-\text{orbit of }(B'_1,B'_2)\m pos(t_Z\i(B'_1),t_Z\i(B'_2))$
\nl
is a map $\t_Z:W'@>>>W$.

Let $B'\in\cb'$ and let $T'$ be a maximal torus of $B'$. Let $B=t_Z\i(B')\in Z$, and let $T$ be the
unique maximal torus of $G^0$ that contains $T'$; we have $T\sub B$. We have a commutative diagram
$$\CD  N_{H^0}(T')/T'@>j>> N_{G^0}(T)/T \\
        @Va'VV                  @VaVV        \\
        W'@>\t_Z>>               W          \\
        \endCD$$
where $a'(nT')=H^0-\text{orbit of } (B',nB'n\i)$, $a(nT)=G^0-\text{orbit of } (B,nBn\i)$ and $j$ is as in
1.5(a). Since $j$ is an injective homomorphism it follows that $\t_Z$ is an injective homomorphism. Now $j$ has
image contained in $(N_{G^0}(T)/T)^h$, see 1.5(b). It follows that $\t_Z$ has image contained in $W^D$
so that $t_Z$ restricts to an (injective) homomorphism $\t_Z^D:W'@>>>W^D$.

Now let $Z,Z'$ in $\un\cb_h$; let $w=w_{Z,Z'}$ be the unique $G^0$-orbit on $\cb\T\cb$ that
contains $(t_Z\i(B'),t_{Z'}\i(B'))$ for any $B'\in\cb'$ (we view $w$ as an element of $W$). As in 1.6 we
see that

(a) $\t_{Z'}(w')w=w\t_Z(w')$ for any $w'\in W'$.
\nl
We show that

(b) $w\in W^D$.
\nl
It is enough to show that if $B'\in\cb'$ then

$(ht_Z\i(B')h\i,ht_{Z'}\i(B')h\i)=(t_Z\i(hB'h\i),t_{Z'}\i(hB'h\i))$
\nl
is in the same $G^0$-orbit as $(t_Z\i(B'),t_{Z'}\i(B'))$. This is clear since $hB'h\i,B'$ are
Borel subgroups of $H^0$ hence are in the same $H^0$-orbit.

From (a),(b) we see that

(c) {\it the injective homomorphism $\t_Z^D:W'@>>>W^D$ defined for $Z\in\un\cb_h$ is
independent of the choice of $Z$, up to composition with an inner automorphism of $W^D$.}

\subhead 1.8 \endsubhead
Assume that $h\in D$ is a unipotent qss element in $G$. Then

(a) {\it $Z_{G^0}(h)$ is connected, see \cite{DM94,1.33}; hence it is equal to $H^0$ where $H=Z_G(h)$;}

(b) {\it if $h'\in D$ is a unipotent qss element in $G$ then $h'$ is $G^0$-conjugate to $h$, see
\cite{Sp82, II, 2.21}.}
\nl
We show:

(c) {\it there exists a pinning of $G^0$ preserved by $\Ad(h)$.}
\nl
If $p=0$ then $h=1$ (so that $D=G^0$) and the result is trivial. Assume now that $p>1$. We fix a pinning of $G^0$ (that is, a choice of a Borel subgroup of $G^0$,
a maximal torus of it and of isomorphisms from $\kk$ to the root subgroups
corresponding to simple roots and their negatives, satisfying some well known
compatibilities); let $\G$  be
the group of all $g\in G$ such that $\Ad(g)$ preserves this pinning (it may permute non-trivially the simple
root subgroups). Clearly $\G$ is a closed subgroup of $G$ which meets any connected component of $G$; let
$\g\in\G\cap D$. The image of $\g_s$ in $G/G^0$ is semisimple and of order dividing $p$ hence $\g_s=1$ and
$\g=\g_u\in\G\cap D$. Now $\Ad(\g_u)$ preserves our pinning hence $\g_u$ is qss in $G$; it is also unipotent
hence by (b) we have $x\g_ux\i=h$ for some $x\in G^0$. Since $\Ad(\g_u)$ preserves our pinning, it follows
that $\Ad(x\g_ux\i)$ preserves some pinning of $G^0$. Hence $\Ad(h)$ preserves some pinning of $G^0$. This
proves (c).

We show:

(d) {\it $W_{H^0}$ is isomorphic to $W^D$.}
\nl
We choose a pinning of $G^0$ preserved by $\Ad(h)$. Let $T\in\ct_{G^0}$ be associated to
this pinning. Let $T'=T\cap Z_G(h)^0$. Recall from 1.5(b) the imbedding $N_{H^0}(T')/T'@>>>(N_{G^0}(T)/T)^h$. 
We show that this is surjective.
Let $\ci$ be the image of the Tits section $N_{G^0}T/T@>>>N_{G^0}T$; this is defined in
terms of the pinning hence $h\ci h\i=\ci$. Let $n\in N_{G^0}T$ be such that $hnh\i\in nT$. We have $n\in n_0T$
where $n_0\in\ci$ and $hn_0h\i\in n_0T\cap\ci$. But if $n_0t\in\ci$, $t\in T$ then $t=1$. Thus $hn_0h\i=n_0$
so that $n_0\in N_{Z_{G^0}(h)}(T)$. Using (a) we deduce $n_0\in N_{H^0}(T)=N_{H^0}(T')$. This shows
that our map is surjective hence bijective. This proves (d).

In our case 1.7(c) implies (using (d)):

(e)  {\it For any $Z\in\un\cb_h$, $\t_Z^D:W_{H^0}@>>>W^D$ is an isomorphism; moreover, it is
independent of the choice of $Z$, up to composition with an inner automorphism of $W^D$.}

\subhead 1.9\endsubhead
To any $g\in D$ we associate $E_g\in\Irr(W^D)$ as follows.

(i) Assume first that $g$ is unipotent. Then $E_g$ is defined as in 1.3(a).

(ii) Next we assume that $g_s$ is central in $G$. Let $D_1$ be the connected component of $G$ that contains
$g_u$. We have $W^D=W^{D_1}$. Then $E_{g_u}\in\Irr(W^{D_1})=\Irr(W^D)$ is defined as in (i) (with $G,D,g$
replaced by $G,D_1,g_u$). We set $E_g=E_{g_u}$. If $g$ is unipotent this definition agrees with the one in (i).

(iii) Consider the general case. Let $D'$ be the connected component of $G':=Z_G(g_s)$ which contains $g$ and
let $W'=W_{G'{}^0}$. Then $W'{}^{D'}$ is defined. Now $g_s$ is central in $G'$. Let $E'_g\in\Irr(W'{}^{D'})$
be defined as $E_g$ in (ii) by replacing $G,D,W,g$ by $G',D',W',g$. Let
$E_g=j_{W'{}^{D'}}^{W^D}(E'_g)\in\Irr(W^D)$ where $j_{-}^{-}(-)$ is the $j$-induction \cite{LS79} from
$W'{}^{D'}$ to $W^D$. (We use any one of the imbeddings of $W'{}^{D'}$ into $W^D$ given by 1.6(c).) Note
that $E'_g$ is a good representation of $W'{}^{D'}$, in the sense of \cite{L15, 0.2} so the $j$-induction can
be applied to it. (To see this we can assume that $G^0$ is simple; if moreover $D=G^0$, we use the explicit
knowledge of $E_g$ in (ii) for connected groups, see \cite{AL82}, \cite{L84a},\cite{LS85},\cite{Sh80},
\cite{Sp85}; if $D\ne G^0$,
all relevant irreducible representations are automatically good.) If $g_s$ is central, this definition agrees
with the one in (ii). (In this case, $G'=G,D'=D,W'=W$.)

\subhead 1.10\endsubhead
In the case where $D_{un}\ne\emp$ we define $\cs_1(G,D)$ to be the subset of $\Irr(W^D)$ consisting of
irreducible representations $E$ such that $E=E_g$ for some $g\in D_{un}$. In the general case,
let $\cs_2(G,D)$ be the subset of $\Irr(W^D)$ consisting of irreducible representations $E$
such that $E=E_g$ for some $g\in D$. Equivalently, $\cs_2(G,D)$ is the subset of $\Irr(W^D)$ consisting of
irreducible
representations $E$ such that $E=E_g$ for some $g\in D$ with $g_s$ of finite order in $G$. (We use the
fact that for any $g\in D$ we can find $g'\in D$ such that
$Z_G(g_s)=Z_G(g'_s)$ and $g'_s$ has finite order in $G$.)
The representations in $\cs_2(G,D)$ are said to be the {\it $2$-special} representations of $W^D$. (This set is
not attached to the Weyl group $W^D$ but rather to the pair $(G,D)$.)
In the case where $D=G^0$, $\cs_2(G,D)$ coincides with the set $\cs_2(W)$ defined in \cite{L15, 1.1}.
This follows from the description of $\cs_2(W)$ given in \cite{L15, 2.1}.

\subhead 1.11\endsubhead
From the definition, the representations of $W^D$ in $\cs_2(G,D)$ are exactly those  obtained by applying
$j_{W_s^{D_1}}^{W^D}$ to any representation of $W_s^{D_1}$ in $\cs_1(Z_G(s),D_1)$ where $s$ is any semisimple
element in $G$ and $D_1$ is any connected component of $Z_G(s)$ such that $(D_1)_{un}\ne\emp$
and $sD_1\sub D$; here $W_s=W_{Z_G(s)^0}$ and $W_s^{D_1}$ is viewed as a subgroup of $W^D$ as in 1.6(c).

\subhead 1.12\endsubhead
We define a partition of $D$ into subsets (called Strata): we say
that $g,g'$ in $D$ are in the same Stratum if $E_g=E_{g'}$. (In the case where $D=G^0$, this partition
coincides with the partition defined in \cite{L15, 2.3}.) Let $Str(D)$ be the set of Strata of $D$.
Now $g\m E_g$ defines a bijection $Str(D)@>\si>>\cs_2(G,D)$. 
One can show that any Stratum of $D$ is a constructible subset of $D$. It is likely that

(b) {\it any Stratum of $D$ is locally closed in $D$.}
\nl
When $D=G^0$ this is proved in \cite{C20}.

\subhead 1.13\endsubhead
Let $\cz_{G^0}$ be the centre of $G^0$. Let $G'=G/\cz_{G^0}$ and let
$D'$ be the image of $D$ under $G@>>>G'$ (a connected component of $G'$). From the definition, the Strata of
$G$ are exactly the inverse images of the Strata of $D'$ under $G@>>>G'$.
Thus $Str(D),Str(D')$ are naturally in bijection. The fixed point set of $[D']$ on
$W_{G'{}^0}=W$ is equal to $W^D$. From the definition we have $\cs_2(G',D')=\cs_2(G,D)$.

\subhead 1.14\endsubhead
Let $(W_1,\g)$ be a pair consisting of a Weyl group $W_1$ and a Coxeter group automorphism $\g:W_1@>>>W_1$
which is ordinary in the sense that, whenever $s\ne s'$ are simple reflections in the same orbit of $\g$, the
product $ss'$ has order $2$ or $3$. Let $W_1^\g$ be the fixed point set of $\g:W_1@>>>W_1$.
We state our main result.

\proclaim{Theorem 1.15} To any $(W_1,\g)$ as in 1.14 one can associate canonically a subset $\cs_2(W_1,\g)$
of $\Irr(W_1^\g)$ such that, whenever $(W_1,\g)=(W,[D])$ with $G,D,W,[D]$ as above, we have
$\cs_2(W_1,\g)=\cs_2(G,D)$.
\endproclaim
The proof is given in \S4.

\subhead 1.16\endsubhead
Following \cite{Sp82, I, 1.1} we set

$\rk_D(G)=\max_{g\in D}(\text{ dimension of a maximal torus of }Z_G(g)^0)$.
\nl
Let $g\in D$. Let $D_1$ be the connected component of $G':=Z_G(g_s)$ which contains $g_u$. We have

(a) $\rk_{D_1}(G')=\rk_D(G)$.
\nl
See \cite{Sp82, II, 1.14}. From \cite{Sp82, II, 10.15} we have

(b) $\dim Z_{G'}(g_u)=2\dim\cb'_{g_u}+\rk_{D_1}(G')$
\nl
where $\cb'_{g_u}$ is the variety of Borel subgroups of $G'{}^0$ which are normalized by $g_u$.
From the theory of Springer correspondence it is known that 
$\dim\cb'_{g_u}=b_{E'_g}$ where $E'_g\in\Irr(W'{}^{D'}),W',D'$ are as in 1.9(iii).
Since the $b$-function is preserved by $j$-induction, we have also $\dim\cb'_{g_u}=b_{E_g}$. Thus (b) becomes
(using also (a)):

(c) $\dim Z_{G'}(g_u)=2b_{E_g}+\rk_D(G)$.
\nl
We have clearly $Z_{G'}(g_u)=Z_G(g)$ and $\dim Z_G(g)=\dim Z_{G^0}(g)$. It follows that

(d) $\dim Z_{G^0}(g)=2b_{E_g}+\rk_D(G)$.
\nl
From this we see that $\dim Z_{G^0}(g)$ is constant when $g$ varies in a fixed Stratum of $D$.
We also see that

(e) {\it any  $\Si\in Str(D)$ is a union of $G^0$-orbits of fixed dimension, namely
$\dim G^0-2b_{E_g}-\rk_D(G)$ where $g\in\Si$.}
\nl
In the case $D=G^0$ this is proved in \cite{L15, 2.4}. 

\head 2. Bipartitions\endhead
\subhead 2.1 \endsubhead
A {\it bipartition} is a sequence $\l=(\l_1,\l_2,\l_3,\do)$ in $\NN$ such that $\l_m=0$ for large $m$ and
$\l_1\ge\l_3\ge\l_5\ge\do$, $\l_2\ge\l_4\ge\l_6\ge\do$. We write $|\l|=\l_1+\l_2+\l_3+\do$. For $n\in\NN$
let $BP^n$ be the set of bipartitions $\l$ such that $|\l|=n$. Let $(e,e')\in\NN\T\NN$. We say that a
bipartition $(\l_1,\l_2,\l_3,\do)$ has {\it excess} $(e,e')$ if $\l_i+e\ge\l_{i+1}$ for $i=1,3,5,\do$ and
$\l_i+e'\ge\l_{i+1}$ for $i=2,4,6,\do$. Let $BP_{e,e'}$ be the set of bipartitions which have excess $(e,e')$.
For $n\in\NN$ let $BP^n_{e,e'}=BP^n\cap BP_{e,e'}$. A {\it partition} is the same as a bipartition of
excess $(0,0)$. A bipartition is the same as an ordered pair of partitions

$(\l_1,\l_3,\l_5,\do),(\l_2,\l_4,\l_6,\do))$.

If $\l=(\l_1,\l_2,\l_3,\do)$,  $\ti\l=(\ti\l_1,\ti\l_2,\ti\l_3,\do)$ are bipartitions, then
$\l+\ti\l:=(\l_1+\ti\l_1,\l_2+\ti\l_2,\l_3+\ti\l_3,\do)$ is a bipartition; moreover if
$(e,e')\in\NN$, $(\te,\te')\in\NN$ and $\l\in BP_{e,e'},\ti\l\in BP_{\te,\te'}$ then
$\l+\ti\l\in BP_{e+\te,e'+\te'}$.

\subhead 2.2 \endsubhead
For $(e,e')\in\NN\T\NN$, $m\ge1$ let ${}^m\hBP_{e,e'}$ be the set of all sequences
$\L=(\L_1,\L_2,\L_3,\do,\L_{2m+1})$ in $\NN$ such that $\L_1\ge\L_2\ge\L_3\ge\do\ge\L_{2m+1}$,
$\L_i-\L_{i+2}\ge e+e'$ for $i=1,2,\do,2m-1$, $\L_{2m}=e'$, $\L_{2m+1}=0$; let ${}^m\hBP'_{e,e'}$ be the set
of all sequences $\L=(\L_1,\L_2,\L_3,\do,\L_{2m})$ in $\NN$ such that $\L_1\ge\L_2\ge\L_3\ge\do\ge\L_{2m}$,
$\L_i-\L_{i+2}\ge e+e'$ for $i=1,2,\do,2m-2$, $\L_{2m-1}=e$, $\L_{2m}=0$.
Given $n\in\NN$ we define ${}^m\k:BP^n_{e,e'}@>>>{}^m\hBP_{e,e'}$ by
$(\l_1,\l_2,\l_3,\do)\m(\L_1,\L_2,\L_3,\do,\L_{2m+1})$ where

$\L_1=\l_1+m(e+e'),\L_2=\l_2+(m-1)e+me',\L_3=\l_3+(m-1)(e+e'),$

$\L_4=\l_4+(m-2)e+(m-1)e',\do,\L_{2m}=\l_{2m}+e',  \L_{2m+1}=\l_{2m+1}$.
\nl
(We choose $m$ large enough so that $\l_{2m}=\l_{2m+1}=0$ for all $\l\in BP^n_{e,e'}$.)
Then ${}^m\k$ defines

(a) a bijection of $BP^n_{e,e'}$ onto a subset ${}^m\hBP^n_{e,e'}$ of ${}^m\hBP_{e,e'}$.

Given $n\in\NN$ we define ${}^m\k':BP^n_{e,e'}@>>>{}^m\hBP'_{e,e'}$ by

$(\l_1,\l_2,\l_3,\do)\m(\L_1,\L_2,\L_3,\do,\L_{2m})$
\nl
where

$\L_1=\l_1+me+(m-1)e',\L_2=\l_2+(m-1)(e+e'),$

$\L_3=\l_3+(m-1)e+(m-2)e',\L_4=\l_4+(m-2)(e+e'),\do,$

$\L_{2m-1}=\l_{2m-1}+e, \L_{2m}=\l_{2m}$.
\nl
(We choose $m$ large enough so that $\l_{2m-1}=\l_{2m}=0$ for all $\l\in BP^n_{e,e'}$.)
Then ${}^m\k'$ defines

(b) a bijection of $BP^n_{e,e'}$ onto a subset ${}^m\hBP'{}^n_{e,e'}$ of ${}^m\hBP'_{e,e'}$.

\proclaim{Proposition 2.3}Let $N\ge2,u\ge2$. Let $c\in[0,N]$. We set $\e=0$ if $c=0$, $\e=1$ if $c>0$.
Let $A_1\le A_2\le\do\le A_u$ be a sequence in $\NN$ such that $A_1=0$, $A_2=c$,
$A_3-A_1\ge N,A_4-A_2\ge N,\do,A_u-A_{u-2}\ge N$. There exist two sequences
$B_1\le B_2\le\do\le B_u$, $C_1\le C_2\le\do\le C_u$ in $\NN$ such that $A_i=B_i+C_i$ for all
$i$, $B_1=0$, $B_2=\e$, $B_3-B_1\ge1,B_4-B_2\ge1,\do,B_u-B_{u-2}\ge1$,
$C_3-C_1\ge N-1,C_4-C_2\ge N-1,...,C_u-C_{u-2}\ge N-1$.
\endproclaim
We say that $i\in[1,u]$ is single if $A_i$ appears exactly once in the sequence $A_1,A_2,...,A_u$.
We define $B_s\in\NN$ by induction on $s$. We set $B_1=0$, $B_2=\e$. Assume that $s\in[3,u]$
and that $B_i$ is already defined for $i\in[1,s-1]$.

(i) If $s$ is single and $s-1$ is not single we set $B_s=B_{s-1}+1$.

(ii) If $s$ and $s-1$ are singles and $s-2$ is not single we set $B_s=B_{s-1}$.

(iii) If $s,s-1,s-2$ are singles and $B_{s-1}=B_{s-2}$ we set $B_s=B_{s-1}+1$.

(iv) If $s,s-1,s-2$ are singles and $B_{s-1}\ne B_{s-2}$ we set $B_s=B_{s-1}$.

(v) If $s$ is not single and $s-1$ is single we set $B_s=B_{s-1}+1$.

(vi) If $s$ is not single, $s-1$ is not single and $A_{s-1}<A_s$ we set $B_s=B_{s-1}+1$.

(vii) If $s$ is not single, $s-1$ is not single and $A_{s-1}=A_s$ we set $B_s=B_{s-1}$.
\nl
This completes the inductive definition of $B_s$.

We show by induction on $s\in[3,u]$ that

(a) $A_s-A_{s-1}-B_s+B_{s-1}\ge0$,

(b) $A_s-A_{s-2}-B_s+B_{s-2}-(N-1)\ge0$.
\nl
Assume that $s$ is as in (i) so that $B_s=B_{s-1}+1=B_{s-2}+1$. We have $A_s\ge A_{s-1}+1$, hence

$A_s-A_{s-1}-B_s+B_{s-1}\ge1-B_s+B_{s-1}\ge0$,

$A_s-A_{s-2}-B_s+B_{s-2}-(N-1)=A_s-A_{s-2}-1-(N-1)\ge N-N=0$.
\nl
Assume that $s$ is as in (ii) so that $B_s=B_{s-1}=B_{s-2}+1$. We have

$A_s-A_{s-1}-B_s+B_{s-1}=A_s-A_{s-1}\ge0$,

$A_s-A_{s-2}-B_s+B_{s-2}-(N-1)=A_s-A_{s-2}-N\ge0$.
\nl
Assume that $s$ is as in (iii) so that $B_s=B_{s-1}+1=B_{s-2}+1$. We have

$A_s\ge A_{s-1}+1$, hence

$A_s-A_{s-1}-B_s+B_{s-1}\ge 1-B_s+B_{s-1}\ge0$,

$A_s-A_{s-2}-B_s+B_{s-2}-(N-1)=A_s-A_{s-2}-N\ge0$.
\nl
Assume that $s$ is as in (iv) so that $B_s=B_{s-1}=B_{s-2}+1$. We have 

$A_s-A_{s-1}-B_s+B_{s-1}=A_s-A_{s-1}\ge0$,

$A_s-A_{s-2}-B_s+B_{s-2}-(N-1)=A_s-A_{s-2}-N\ge0$.
\nl
Assume that $s$ is as in (v) so that $B_s=B_{s-1}+1$.

Since $s-1$ is single we have $A_s\ge A_{s-1}+1$, hence

$A_s-A_{s-1}-B_s+B_{s-1}\ge 1-B_s+B_{s-1}=0$.

Since $s$ is not single and $s-1$ is single we have

$A_{s+1}=A_s>A_{s-1}>A_{s-2}$.
\nl
We have $B_{s-1}=B_{s-2}$ or $B_{s-1}=B_{s-2}+1$. If $B_{s-1}=B_{s-2}$ then

$A_s-A_{s-2}-B_s+B_{s-2}-(N-1)\ge N-B_{s-1}-1+B_{s-1}-(N-1)=0$.
\nl
If $B_{s-1}=B_{s-2}+1$ then 

$A_s-A_{s-2}-B_s+B_{s-2}-(N-1)=A_{s+1}-A_{s-2}-2-(N-1)
=(A_{s+1}-A_{s-1})+(A_{s-1}-A_{s-2})-N-1\ge N+1-N-1=0$.
\nl
Assume that $s$ is as in (vi) so that $B_s=B_{s-1}+1$, $B_{s-1}=B_{s-2}$, $A_s>A_{s-1}$. We have

$A_s-A_{s-1}-B_s+B_{s-1}=A_s-A_{s-1}-1\ge0$,

$A_s-A_{s-2}-B_s+B_{s-2}-(N-1)=A_s-A_{s-2}-1-(N-1)\ge0$.
\nl
Assume that $s$ is as in (vii) so that $B_s=B_{s-1}$, $B_{s-1}=B_{s-2}+1$. We have

$A_s-A_{s-1}-B_s+B_{s-1}=0$,

$A_s-A_{s-2}-B_s+B_{s-2}-(N-1)=A_s-A_{s-2}-1-(N-1)\ge0$.
\nl
This completes the inductive proof of (a),(b).
We set $C_i=A_i-B_i$. Then $B_i,C_i$ satisfy the requirements of the proposition.
(We have $C_1=0$, $C_2=c-\e\ge C_1$.) The proposition is proved. 

\proclaim{Corollary 2.4} Let $(e,e')\in\NN\T\NN$, $(\te,\te')\in\NN\T\NN$. Let $\l\in BP_{e+\te,e'+\te'}$.
There exist $\l^1\in BP_{e,e'}$, $\ti\l^1\in BP_{\te,\te'}$ such that $\l=\l^1+\ti\l^1$.
\endproclaim
Using 2.2(a) with $m\gg0$, we see that it is enough to show that

(a) for any $\L\in{}^m\hBP_{e+\te,e+\te'}$ there exist $\L^1\in{}^m\hBP_{e,e'}$,
$\ti\L^1\in{}^m\hBP_{\te,\te'}$ such that $\L=\L^1+\ti\L^1$ (addition coordinate by coordinate).
\nl
The following statement clearly implies (a):

(b) for any $(f,f')\in\NN\T\NN$ and any $\L=(\L_1,\L_2,\do,\L_{2m+1})\in{}^m\hBP_{f,f'}$ there exist
$\L^1,\L^2,\do,\L^f$ in ${}^m\hBP_{1,0}$ and $\L'{}^1,\L'{}^2,\do,\L'{}^{f'}$ in ${}^m\hBP_{0,1}$ and
$\ti\L^0\in{}^m\hBP_{0,0}$ such that $\L=\L^1+\L^2+\do+\L^f+\L'{}^1+\L'{}^2+\do+\L'{}^{f'}+\ti\L^0$.
\nl
We prove (b). Using 2.3 with $u=2m+1$, $N=f+f'$ and $(A_1,A_2,\do,A_{2m+1})=(\L_{2m+1},\L_{2m},\do,\L_1)$,
$c=f'$, we see that if $f'\ge1$, then $\L=\L'{}^1+\ti\L$ where $\L'{}^1\in{}^m\hBP_{0,1}$,
$\ti\L\in\hBP_{f,f'-1}$.
Using this repeatedly, we see that $\L$ is of the form $\L'{}^1+\L'{}^2+\do+\L'{}^{f'}+\ti\L^1$ where
$\L'{}^1,\L'{}^2,\do,\L'{}^{f'}$ are in ${}^m\hBP_{0,1}$ and
$\ti\L^1\in\hBP_{f,0}$. Thus it is enough to prove (b) assuming in addition that $f'=0$.

Using 2.3 with $u=2m+1$, $N=f$ and

$(A_1,A_2,\do,A_{2m+1})=(\L_{2m+1},\L_{2m},\do,\L_1)$,
\nl

$c=0$, we see that if $f\ge1$, then $\L=\L^1+\ti\L$ where $\L^1\in{}^m\hBP_{1,0}$, $\ti\L\in\hBP_{f-1,0}$.
Using this repeatedly, we see that $\L$ is of the form
$\L^1+\L^2+\do+\L^f+\ti\L^0$ where $\L^1,\L^2,\do,\L^f$ are in ${}^m\hBP_{1,0}$ and
$\ti\L^0\in\hBP_{0,0}$. This proves (b). The corollary is proved.

\subhead 2.5\endsubhead
Some special cases of 2.4 have been used (without proof) in three surjectivity statements in \cite{L15}
(see \cite{L15, p.348, line 8},\cite{L15, p.349, line 16 of 3.7}, \cite{L15, p.351, line 5}), namely that
addition defines surjective maps $BP_{1,1}\T BP_{1,1}@>>>BP_{2,2}$, $BP_{2,0}\T BP_{0,2}@>>>BP_{2,2}$,
$BP_{0,2}\T BP_{0,2}@>>>BP_{0,4}$. Other special cases of 2.4 are contained in \cite{L09}, namely that
addition defines surjective maps $BP_{1,0}\T BP_{0,1}@>>>BP_{1,1}$, $BP_{1,0}\T BP_{1,0}@>>>BP_{2,0}$,
$BP_{0,1}\T BP_{0,1}@>>>BP_{0,2}$.

\head 3. Examples\endhead
\subhead 3.1\endsubhead
Assume that $p\ne2$, that $G=O(V)$ where $V$ is a $\kk$-vector space of even dimension $N\ge4$
with a fixed nondegenerate symmetric bilinear form and that $D=O(V)-SO(V)$. Let $g\in D$. For any $c\in\kk^*$
let $V_c$ be the $c$-eigenspace of $g_s:V@>>>V$. Let $d_c=\dim V_c$. For any $c\in\kk^*$
such that $c^2\ne1$ let $\l^c_1\ge\l^c_2\ge\do$ be the partition of $d_c$ whose nonzero terms are the sizes of
the Jordan blocks of $g_u:V_c@>>>V_c$. For $c\in\kk^*$ such that $c^2=1$ let $\nu^c_1\ge\nu^c_2\ge\do$ be
the partition of (the odd number)
$d_c$ whose nonzero terms are the sizes of the Jordan blocks of the unipotent element
$g_u\in SO(V_c)$. Let $\l^c=(\l^c_1,\l^c_2,\l^c_3,\do)$ be
the bipartition in $BP^{(d_c-1)/2}_{2,0}$ associated to $\nu^c_1\ge\nu^c_2\ge\do$ in \cite{L15, 3.6(a)}.
Note that $\l^c$ is such that the Springer representation attached to the unipotent element
$g_u\in SO(V_c)$ (an irreducible representation of a Weyl group of type $B$)
is indexed by $\l^c$. (We use results in \cite{L84a}.) Define ${}^g\l=({}^g\l_1,{}^g\l_2,{}^g\l_3,\do)$ by
${}^g\l_j=\sum_c \l^c_j$ where $c$ runs over a set of representatives for the orbits of the involution
$a\m a\i$ of $\kk^*$. Note that ${}^g\l\in BP^{(N-2)/2}_{4,0}$. Thus $g\m{}^g\l$  defines a map
$D@>>>BP^{(N-2)/2}_{4,0}$. From the definitions we see that the fibres of this map are precisely the Strata
of $D$. (We use the description of the $j$-induction in classical Weyl groups groups given in \cite{L09}.)
This map is also surjective: indeed, by 2.4, for any $\mu\in BP^{(N-2)/2}_{4,0}$ we can find
$\mu'\in BP_{2,0},\mu''\in BP_{2,0}$ such that $\mu=\mu'+\mu''$; we have $\mu'\in BP^k_{2,0}$,
$\mu''\in BP^{k'}_{2,0}$ for some $k\ge0,k'\ge0$ such that $k+k'=(N-2)/2$ and it remains to use that any
bipartition in $BP^k_{2,0}$ (resp. $BP^{k'}_{2,0}$) represents the Springer representation associated to a
unipotent element in $SO_{2k+1}(\kk)$ (resp. $SO_{2k'+1}(\kk)$); see \cite{L15, 3.6(a)}. We see that

(a) $D\m BP^{(N-2)/2}_{4,0}$, $g\m{}^g\l$ defines a bijection from the set of Strata of $D$ to
$BP^{(N-2)/2}_{4,0}$. Thus, $\cs_2(G,D)$ can be identified with $BP^{(N-2)/2}_{4,0}$. 

\subhead 3.2\endsubhead
Assume that $V$ is a $\kk$-vector space of dimension $N$. Let $V^*$ be
the vector space dual to $V$. For $x\in V,\x\in V^*$ we set $(x,\x)=\x(x)$. Let $D$
be the set of all linear isomorphisms $g:V@>\si>>V^*$. For $g\in D$ we define an isomorphism
$\hg:V^*@>>>V$ by $(\hg(\x),g(x))=(x,\x)$ for any $x\in V,\x\in V^*$. For $g\in GL(V)$ we
define $\hg\in GL(V^*)$ by $(g(x),\hg(\x))=(x,\x)$ for any $x\in V,\x\in V^*$. Let
$G=GL(V)\sqc D$. For $g,g'$ in $G$ we define $g*g'\in G$ to be $gg'$ if $g'\in GL(V)$ and
$\hg g'$ if $g'\in D$. Then $(g,g')\m g*g'$ defines a group structure on $G$. Note that
$G$ has an obvious structure of algebraic group over $\kk$ with $G^0=GL(V)$
and with $D$ being another connected component. We now assume that this $G,D$ is the same as that in 1.1.

\subhead 3.3\endsubhead
In the setup of 3.2 assume that $p\ne2$. Let $s\in D$ be semisimple in $G$. Define a bilinear form
$<,>_s:V\T V@>>>\kk$ by $<x,x'>_s=(x,s(x'))$. For any $c\in\kk^*$ let $V_c$ be the
$c$-eigenspace of $s*s:V@>>>V$; we have

$V_c=\{x\in V;<y,x>_s=c<x,y>_s\text{ for any }y\in V\}$;
\nl
let $d_c=\dim V_c$. It follows that if $c\in\kk^*,c'\in\kk^*,cc'\ne1$ and $x\in V_c,y\in V_{c'}$ then
$<x,y>_s=cc'<x,y>_c$ so that $<x,y>_s=0$. It follows that $<,>_s$ defines an isomorphism
$V_c@>>>V_{c\i}^*$. If $h\in Z_{G^0}(s)$ then $h$ restricts to an isomorphism $h_c:V_c@>>>V_c$ for any
$c\in\kk^*$. Now $h_{-1}:V_{-1}@>>>V_{-1}$ preserves the nondegenerate symplectic form $<,>_s$
restricted to $V_{-1}$; $h_1:V_1@>>>V_1$ preserves the nondegenerate symmetric bilinear form $<,>_s$
restricted to $V_1$; if $c\in\kk^*-\{1,-1\}$ then $h_c:V_c@>>>V_c$ is related to $h_{c\i}:V_{c\i}@>>>V_{c\i}$
by $<h_c(x),h_{c\i}(y)>_s=<x,y>_s$ for $x\in V_c,y\in V_{c\i}$; moreover, $h\m(h_c)$ is an isomorphism of
$Z_{G^0}(s)$ onto $Sp(V_{-1})\T O(V_1)\T\prod_{c\in C'}GL(V_{c'})$ where $C'$ is a set of representatives
for the orbits of $c\m c\i$ on $\kk^*-\{1,-1\}$.

Now assume that $g\in D$ and $s=g_s$. Then $g_u\in Z_{G^0}(s)$ gives rise to

a unipotent element $(g_u)_{-1}\in Sp(V_{-1})$ with Jordan blocks of sizes
$\nu_1\ge\nu_2\ge\do$ (a partition of $d_{-1}$);

a unipotent element $(g_u)_1\in O(V_1)$ with Jordan blocks of sizes
$\nu'_1\ge\nu'_2\ge\do$ (a partition of $d_1$);

a unipotent element $(g_u)_c\in GL(V_c)$ with Jordan blocks of sizes
$\l^c_1\ge\l^c_2\ge\do$ (a partition of $d_c$) for any $c\in C'$.
\nl
Let $(\l_1,\l_2,\l_3,\do)$ be the bipartition in $BP^{d_{-1}/2}_{1,1}$ associated to
$\nu_1\ge\nu_2\ge\do$ in \cite{L15, 3.4(a)}. Note that $(\l_1,\l_2,\l_3,\do)$ is such that the Springer
representation attached to the unipotent element $g_u\in Sp(V_{-1})$ (an irreducible representation of a
Weyl group of type $C$) is indexed by $(\l_1,\l_2,\l_3,\do)$. (We use results in \cite{L84a}.)

Let $(\l'_1,\l'_2,\l'_3,\do)$ be the bipartition in $BP^{(d_1-1)/2}_{2,0}$, if $d_1$ is odd (resp. in
$BP^{d_1/2}_{0,2}$, if $d_1$ is even) associated to $\nu'_1\ge\nu'_2\ge\do$ in \cite{L15, 3.6(a)} (resp.
\cite{L15. 3.6(b)}). Note that $(\l'_1,\l'_2,\l'_3,\do)$ is such that the Springer representation attached
to the unipotent element $g_u\in SO(V_1)$ (an irreducible representation of a Weyl group of type $B$ or $D$)
is indexed by $(\l'_1,\l'_2,\l'_3,\do)$. (We use results in \cite{L84a}.)
Define ${}^g\l=({}^g\l_1,{}^g\l_2,{}^g\l_3,\do)$ by ${}^g\l_j=\l_j+\l'_j+\sum_{c\in C'}\l^c_j$.
We have ${}^g\l\in BP^{(N-1)/2}_{3,1}$ (if $N$ is odd), ${}^g\l\in BP^{N/2}_{1,3}$ (if $N$ is even). Thus,
$g\m{}^g\l$ defines a map $D@>>>BP^{(N-1)/2}_{3,1}$ (if $N$ is odd) or $D@>>>BP^{N/2}_{1,3}$ (if $N$ is
even). From the definitions we see that the fibres of this map are precisely the Strata of $D$. (We use the
description of the $j$-induction in classical Weyl groups groups given in \cite{L09}.)

This map is also surjective: indeed, by 2.4, for any $\mu\in BP^{(N-1)/2}_{3,1}$ if $N$ is odd
(resp. $\mu\in BP^{N/2}_{1,3}$ if $N$ is even) we can find $\mu'\in BP_{1,1},\mu''\in BP_{2,0}$ if $N$ is odd
(resp. $\mu'\in BP_{1,1},\mu''\in BP_{0,2}$ if $N$ is even) such that $\mu=\mu'+\mu''$; we have
$\mu'\in BP^k_{1,1}$ and $\mu''\in BP^{k'}_{2,0}$ if $N$ is odd (resp. $\mu''\in BP^{k'}_{0,2}$ if $N$ is
even) for some $k\ge0,k'\ge0$ such that $k+k'=(N-1)/2$ if $N$ is odd (resp. $k+k'=N/2$ if $N$ is even)
and it remains to use that any bipartition in $BP^k_{2,0}$ (resp. $BP^{k'}_{2,0}$ if $N$ is odd or
$BP^{k'}_{0,2}$ if $N$ is even) represents the Springer representation associated to a
unipotent element in $Sp_{2k}(\kk)$ (resp. $SO_{2k'+1}(\kk)$ if $N$ is odd or $SO_{2k'}(\kk)$ if $N$ is even);
see \cite{L15, 3.4(a)},\cite{L15, 3.6(a)},\cite{L15, 3.6(b)}. We see that

(a) the map $g\m{}^g\l$ from $D$ to $BP^{(N-1)/2}_{3,1}$ (if $N$ is odd) or to $BP^{N/2}_{1,3}$ (if $N$ is
even) defines a bijection from the set of Strata of $D$ to $BP^{(N-1)/2}_{3,1}$ (if $N$ is odd) or to
$BP^{N/2}_{1,3}$ (if $N$ is even). Thus, $\cs_2(G,D)$ can be identified with $BP^{(N-1)/2}_{3,1}$ (if $N$ is
odd) or to $BP^{N/2}_{1,3}$ (if $N$ is even). 

\subhead 3.4\endsubhead
In the setup of 3.2, with $N=3$,
$D$ is a union of two Strata, $\Si$ and $\Si'$. Now $\Si$ is the union of all
$G^0$-conjugacy classes of dimension $8$ in $D$; $\Si'$ is the union of all $G^0$-conjugacy classes of
dimension $6$ in $D$. More precisely, if $p\ne2$, $\Si'$ consists of two semisimple $G^0$-orbits in $D$; if
$p=2$, $\Si'$ consists of a single $G^0$-orbit in $D$ (its elements are unipotent, qss). We have
$\cs_2(G,D)=\Irr(W^D)$.

\head 4. Proof of Theorem 1.15\endhead
\subhead 4.1\endsubhead
In this section we assume (until the end of 4.11)
that we are given $\vt\in G$ of finite order $\d$ such that
$G=\sqc_{i=1}^\d G^0\vt^i$, that $D=G^0\vt$, that $G^0$ is almost simple, simply connected, with a fixed
pinning, that $\Ad(\vt):G^0@>>>G^0$ has order $\d$ and it preserves the pinning of $G^0$.
(We must have $\d\in\{1,2,3\}$.)

We now assume (until  the end of 4.11) that $\d\ge2$. We are in one of the following cases
(in each case we indicate the type of $G^0$; the number $\d$ appears as an upper index): type
${}^2A_{2n+1},n\ge1$; type ${}^2A_{2n},n\ge1$; type ${}^2D_n,n\ge4$; type ${}^3D_4$; type ${}^2E_6$.

Let $(T,B)\in\tcb_\vt$. Let $\cy=\Hom(\kk^*,T),\cx=\Hom(T,\kk^*)=\Hom(\cy,\ZZ)$. Let $\ucy=\QQ\ot\cy$,
$\ucx=\QQ\ot\cx=\Hom_\QQ(\ucy,\QQ)$. Let $\che R\sub\cy\sub\ucy$ (resp. $R\sub\cx\sub\ucx$) be the set of coroots
(resp. roots) of $G^0$ with respect to $T$. The natural bijection $\che R\lra R$ is denoted by $\che\a\lra\a$.
Let $\{\a_i;i\in I\}$ be the set of simple roots in $R$ determined by $B$. (The corresponding root subgroups are
contained in $B$.) Now $Ad(\vt)$ acts naturally on $\cy$ and $\cx$; this induces a permutation of $R$, one of
$\che R$, one of $\{\a_i;i\in I\}$, and one of $I$ (these permutations are denoted by $\vt$). Let
$G^\vt=Z_G(\vt)$, a reductive group with identity component $G^{\vt0}$. Let $T'=T\cap G^{\vt0}$; we have
$T'\in\ct_{G^{\vt0}}$. We have $T=Z_{G^0}(T')$. Let ${}'\cy=\Hom(\kk^*,T')\sub\cy$,
${}'\cx=\Hom(T',\kk^*)=\Hom({}'\cy,\ZZ)$. Let ${}'\ucy=\QQ\ot{}'\cy$, ${}'\ucx=\QQ\ot{}'\cx$. The homomorphism
$res:\cx@>>>{}'\cx$ (restriction to ${}'\cy$) induces a linear surjective map $res:\ucx@>>>{}'\ucx$. Let
${}'R=res(R)\sub{}'\cx\sub{}'\ucx$. Now $res:R@>>>{}'R$ induces a bijection from the set of orbits of
$\vt:R\to R$ to ${}'R$. For $\b\in{}'R$ let $d'_\b$ be the cardinal of the corresponding $\vt$-orbit on $R$; we
set $d''_\b=2$ if either $2\b\in{}'R$ or $(1/2)\b\in{}'R$ and $d''_\b=1$, otherwise; let $d_\b=d'_\b d''_\b$. For
$\a\in R,\b=res(\a)$, we set:

${}'\che\b=\che\a$ if $d'_\b=1$ or if $d''_\b=2$,

${}'\che\b=\che\a+\vt(\che\a)$ if $d'_\b=2$ and $d''_\b=1$,

${}'\che\b=\che\a+\vt(\che\a)+\vt^2(\che\a)$ if $d_\b=3$,

${}'\che\b=2\che\a+2\vt(\che\a)$ if $d_\b=4$.
\nl
Let ${}'\che R=\{{}'\che\b;\b\in{}'R\}\sub{}'\cy\sub{}'\ucy$; this set is in obvious bijection with
${}'R$ and
$({}'\che R,{}'R,{}'\ucy,{}'\ucx)$ is a (not necessarily reduced) root system. Let ${}'R_0$ be the set of
elements in ${}'R$ which are not in $2{}'R$; let ${}'\che R_0$ be the set of elements in ${}'\che R$ which
are not in $(1/2){}'\che R$. Then $({}'\che R_0,{}'R_0,{}'\ucy,{}'\ucx)$ is a (reduced) root system.
Let $\bI$ be the set of orbits of the bijection $\vt:I\to I$. For $i\in\bI$ let $\b_i=res(\a_{i'})$ where
$i'$ is any element of the orbit $i$. Then $\{\b_i;i\in\bI\}$ is a basis of the root system
$({}'\che R_0,{}'R_0,{}'\ucy,{}'\ucx)$. Let $\uR$ (resp. $\che{\uR}$) be the subset of ${}'\ucx$ (resp. ${}'\ucy$)
consisting of the vectors $d_\b\b$ (resp. $d_\b\i\che\b$) for various $\b\in{}'R$. Then
$(\che{\uR},\uR,{}'\ucy,{}'\ucx)$ is a (reduced, irreducible) root system with basis
$\{\g_i=d_{\b_i}\b_i;i\in\bI\}$. There is a unique vector space isomorphism ${}'\ucx@>>>{}'\ucy$ which
carries
$\g_i$ to ${}'\che\b_i$ for any $i\in\bI$; from the definitions we see that this isomorphism carries
$\uR$ onto ${}'\che R_0$. Hence it carries $\g_0\in\uR$, the negative of the highest root in $\uR$ relative
to the basis $\{\g_i=d_{\b_i}\b_i;i\in\bI\}$, to $h_0$, the negative of the highest coroot in ${}'\che R_0$
relative to the basis $\{{}'\che\b_i;i\in\bI\}$. We have $h_0=\che\b_0$ and $\g_0=d_{\b_0}\b_0$ for a
well defined $\b_0\in{}'R_0$. Setting $\tI=\bI\sqc\{0\}$ there are unique integers $n_i\in\ZZ_{>0}$,
$i\in\tI$ such that $n_0=1$ and $\sum_{i\in\tI}n_i\g_i=0$.
We define a subgroup $\QQ_*$ of $\QQ$ as follows: if $p=0$ then $\QQ_*=\QQ$; if $p>1$, then $\QQ_*$
consists of the rational numbers $q$ such that $Nq\in\ZZ$ for some integer $N$ not divisible by $p$.
Let ${}'\ucy_*=\QQ_*\ot{}'\cy\sub{}'\ucy$. We define
$$\cc=\{y\in{}'\ucy_*;\sum_{i\in\tI}n_i\g_i(y)=0,\g_i(y)\in\QQ_{\ge0}\text{ for }i\in\bI,\g_0(y)+1\in
\QQ_{\ge0}\}.$$
For $y\in\cc$ let $\supp(y)\sub\tI$ be the set of all $i\in\tI$ such that either $i\in\bI,\g_i(y)>0$ or
$i=0,\g_0(y)>-1$. We have $\cc=\sqc_{K\in\cp(\tI)}\cc_K$ where $\cp(\tI)$ is the set of all subsets
$K\sub\tI$ such that $K\ne\emp$ and $\cc_K=\{y\in\cc;\supp(y)=K\}$. 

We define a subset $\cp'(\tI)$ of $\cp(\tI)$ as follows. If $p=0$ we have $\cp'(\tI)=\cp(\tI)$. If $p>1$,
$\cp'(\tI)$ is the set of all $K\in\cp(\tI)$ such that $n_i/p\n\ZZ$ for some $i\in K$.

\proclaim{Lemma 4.2} Assume that $p>1$ and $K\in\cp(\tI)-\cp'(\tI)$. Then $\cc_K=\emp$.
\endproclaim
For any $y\in{}'\ucy_*$ and any $j\in\bI$ we have $\g_j(y)\in\QQ_*$ since $\g_j$ takes integer values on
$\cy$. Assume now that $y\in\cc_K$. Since $n_0=1$, we have $0\n K$ so that $\g_0(y)=-1$. For $i\in K$ we have
$\g_i(y)=s_i/t_i$, $n_i=pn'_i$  where $s_i,t_i,n'_i$ are integers $>0$ with $p$ not dividing $t_i$. We have
$\sum_{i\in K}pn'_is_i/t_i-1=0$. It follows that $p$ divides $\prod_{i\in K}t_i$ hence for some $i\in K$, $p$
divides $t_i$, contradiction.

\proclaim{Lemma 4.3} Assume that $\d\ne p$ and that $K\in\cp'(\tI)$. Then $\cc_K\ne\emp$.
\endproclaim
Assume first that $0\n K$. If $p=0$ we choose any $i\in K\cap\bI$. If $p>1$ we choose $i\in K\cap\bI$ such that
$n_i/p\n\ZZ$. Let $K'=K-\{i\}$. We define $c_j\in\QQ_*$ for $j\in\tI$ as follows: $c_0=-1$; if
$j\in\tI-K,j\ne0$ then $c_j=0$; if $j\in K'$, $c_j\in\QQ_*$ is chosen so that $0<c_j<1/N$ where $N$ is an
integer such that $\sum_{j'\in K'}n_{j'}<N$; $c_i=(1-(\sum_{j'\in K'}n_{j'}c_{j'})/n_i$. (Note that
$c_i\in\QQ_*$ and $c_i>(1-(\sum_{j'\in K'}n_{j'}/N)/n_i>0$.)

Next we assume that $0\in K$. We define $c_j\in\QQ_*$ for $j\in\tI$ as follows: if $j\in\tI-K$ then $c_j=0$;
if $j\in K\cap\bI$ then $c_j\in\QQ_*$ is chosen so that $0<c_j<1/N$ where $N$ is an integer such that
$\sum_{j'\in K-\{0\}}n_{j'}<N$; $c_0=-\sum_{j'\in K\cap\bI}n_{j'}c_{j'}$. (Note that $c_0\in\QQ_*$ and
$c_0>-\sum_{j'\in K\cap\bI}n_{j'}/N>-1$.)

We now define $y\in{}'\ucy$ by $\g_h(y)=c_h$ for all $h\in\bI$. Let ${}'\cx_1$ be the subgroup of ${}'\ucx$
generated by $\{\g_i;i\in\bI\}$. Since $c_h\in\QQ_*$, we see that for any $\g\in{}'\cx_1$ we have
$\g(y)\in\QQ_*$. Using a case by case check we see that the index of ${}'\cx$ in ${}'\cx_1$ is of the form
$\d^m$ for some integer $m\ge1$. It follows that for any $\z\in{}'\cx$ we have $\d^m\z\in{}'\cx_1$ hence
$\d^m\z(y)\in\QQ_*$ hence $\z(y)\in\QQ_*$. (Here we use that $\d\ne p$.) Since this holds for any
$\z\in{}'\cx$ it follows that $y\in{}'\ucy_*$ so that $y\in\cc_K$. This completes the proof.

\subhead 4.4\endsubhead
Let $\mu(\kk)$ be the group of roots of $1$ in $\kk$. We assume that an identification $\QQ_*/\ZZ$ with
$\mu(\kk)$ (as groups) is fixed. Then ${}'\ucy_*/{}'\cy=(\QQ_*/\ZZ)\ot{}'\cy$ is identified with
$\mu(\kk)\ot{}'\cy$ which can be identified (via $\l\ot y\m y(\l)$) with the group $T'_{fin}$ of elements of
finite order in $T'$; thus we obtain a homomorphism $\io:{}'\ucy_*@>>>T'$ whose kernel is ${}'\cy$ and whose
image is $T'_{fin}$. The following result is an adaptation of results in \cite{L02,\S6}. (In {\it loc.cit} the
characteristic is assumed to be zero, but the same arguments apply assuming only that $\d\ne p$, by replacing
the Lie algebras of $T$ and $T'$ by $\ucy,{}'\ucy$.) Results of this kind go back to de Siebenthal's
paper \cite{dS56} where conjugacy classes in disconnected compact Lie groups are discussed.

(a) {\it Assume that $\d\ne p$. Then $x\m\io(x)\vt$ defines a bijection of $\cc$ with a set of representatives
for the $G^0$-conjugacy classes of semisimple elements in $D$ which have finite order in $G$.}
\nl
Let $K\in\cp'(\tI),J=\tI-K$. Let $\che R_J$ be the set of all vectors in ${}'\che R$ which are $\ZZ$-linear
combinations of vectors in $\{\che\b_i;i\in J\}$. Let $R_J$ be the set of all $\b\in{}'R$ such that
${}'\che\b\in\che R_J$. Then for any $x\in\cc_K$, $G_J:=Z_{G^0}(\io(x)\vt)$ is a connected
reductive group (see \cite{St68}) which depends only on $J$, not on $x$,
and $(\che R_J,R_J,{}'\ucy,{}'\ucx)$ is the root system of $G_J$ relative to the maximal torus $T'$ of $G_J$.
(Here $R_J$ is the set of roots and $\che R_J$ is the set of coroots.)

\proclaim{Lemma 4.5}Assume that $\d=p$. Let $s\in G^{\vt0}$ be semisimple. Then $W_{Z_G(s\vt)^0}$
is isomorphic to the fixed point set of $\Ad(\vt)$ on $W_{Z_G(s)^0}$.
\endproclaim
Let $B,T$ in $G^0$ be preserved by $\Ad(\vt)$. Then $(T\cap G^{\vt0},B\cap G^{\vt0})\in\tcb_{G^{\vt0}}$ and
are preserved by $\Ad(\vt)$. Hence $\vt$ is qss in $Z_G(s)$; it is also unipotent. We now use 1.8(a),(d)
with
$G$ replaced by $Z_G(s)$ and note that $Z_G(s)^0\cap Z_G(\vt)=(Z_G(s)^0\cap Z_G(\vt))^0=Z_G(s\vt)^0$.

\proclaim{Lemma 4.6}Assume that $\d=p$. Let $g\in D$. Then some $G^0$-conjugate of $g_s$ is in $Z_G(\vt)^0$.
\endproclaim
The image of $g_s$ in $G/G^0$ is semisimple; since $G/G^0$ is a unipotent group this image must be $1$; thus,
we have $g_s\in G^0$. Let $D'\sub G'=Z_G(g_s)$ be as in 1.6. Now
$G'/G'{}^0@>>>G/G^0$ is injective (and
carries $D'$ to $D$) since $Z_{G^0}(g_s)$ is connected. (Recall that $G^0$ is simply connected.) It follows
that $D'$ has order $p$ in $G'/G'{}^0$. Hence we can find $u\in D'$ unipotent and $(T',B')\in(\tcb_{G'{}^0})_u$.
Let $Z$ be an irreducible component of $\cb_{g_s}$ such that $g(Z)=Z$. Let $B\in Z$ be the unique Borel
subgroup such that $B\cap G'{}^0=B'$. Now $uBu\i\in Z$ (since $g(Z)=Z$) and $uBu\i\cap G'{}^0=uB'u\i=B'$. By
uniqueness we have $uBu\i=B$. Let $T$ be a maximal torus of $B$ containing $T'$. Now there is a unique maximal
torus of $G^0$ containing $T'$. It must be equal to $T$ and $uTu\i$ is a maximal torus containig $T'$ hence
is equal to $T$. Thus $uTu\i=T$. We have $u\in D$ and $u$ is unipotent and qss in $G$. Since $u\in D'$,
we have $g_su=ug_s$. By 1.8(b), for some $h\in G^0$ we have
$huh\i=\vt$, Let $s'=hg_sh\i$. Then $s'\vt=\vt s'$.
Thus we have $s'\in G^0\cap Z_G(\vt)=Z_{G^0}(\vt)=Z_G(\vt)^0$ (we use 1.8(a)).

\subhead 4.7\endsubhead
Assume that $\d\ne p$. Let $g\in D$ be an element of finite order.
By replacing $g$ by a $G^0$-conjugate we can assume that for some
$K\in\cp'(\tI)$ we have $g_s=\io(x)\vt$ where $x\in\cc_K$ (see 4.1, 4.4). Let $H=Z_G(g_s)$. With notation of
4.4 we have $Z_{G^0}(g_s)=G_{\tI-K}=G_{\tI-K}^0=H^0$. Let $W_{\tI-K}=W_{G_{\tI-K}}$. We have $g_s\in D$ hence
$g_u\in G^0$. Thus $g_u\in G^0\cap H$ so that $g_u\in H^0=G_{\tI-K}$. From 1.11 we now see that

(a) $\cs_2(G,D)=\cup_{K\in\cp'(\tI)}j_{W_{\tI-K}}^{W^D}(\cs_1(G_{\tI-K},G_{\tI-K}))$
\nl
where $W_{\tI-K}$ is viewed as a subgroup of $W^D$ as in 1.6(c).
Now let $K\in\cp'(\tI)$. We can find $i\in K$ such that $\{i\}\in\cp'(\tI)$. (If $p=0$ any $i\in K$
satisfies our requirement; if $p>0$ at least one $i\in K$ satisfies our requirement.) We have
$G_{\tI-K}\sub G_{\tI-\{i\}}$ and in fact $G_{\tI-K}$ is a Levi subgroup of a parabolic subgroup
of $G_{\tI-\{i\}}$ so that we can regard $W_{\tI-K}$ as a subgroup of $W_{\tI-\{i\}}$.
If $E\in\cs_1(G_{\tI-K},G_{\tI-K})$ then $j_{W_{\tI-K}}^{W^D}(E)=j_{W_{\tI-\{i\}}}^{W^D}(E')$
where $E'=j_{W_{\tI-K}}^{W_{\tI-\{i\}}}(E)$. It is known that $E'\in\cs_1(G_{\tI-K},G_{\tI-K}))$
(a property of induced unipotent classes). It follows that the union in (a) can be restricted to
one elements subsets $K$ in $\cp'(\tI)$. Thus, we have

(b) $\cs_2(G,D)=\cup_{i\in\tI;\{i\}\in\cp'(\tI)}j_{W_{\tI-\{i\}}}^{W^D}(\cs_1(G_{\tI-\{i\}},G_{\tI-\{i\}}))$.

\subhead 4.8\endsubhead
Assume that $\d=p$. Let $X$ be a set of representatives for the semisimple $G^{\vt0}$-conjugacy classes in
$G^{\vt0}$. Let $g\in D$. By replacing $g$ by a $G^0$-conjugate we can assume that
$g_s\in X$ (see Lemma 4.6). We set $s=g_s$, $H=Z_G(s)$. Since $s\in G^0$ we have $g_u\in D$.
Let $D_s$ be the connected component of $H$ which contains $\vt$; we have $g_u\in D_s$ (indeed,
$g_u\vt\i\in G^0$ and $g_u\vt\i\in H$ hence $g_u\vt\i\in G^0\cap H$ which equals $H^0$ by \cite{St68}). Hence
$W_{H^0}^{D_s}$ is equal to the fixed point of $\Ad(\vt)$ on $W_{H^0}$ which by Lemma 4.5 is the same as
$W_{Z_G(s\vt)^0}$ that is, $W_{Z_{G^{\vt0}}(s)^0}$. (Indeed,
$Z_G(s\vt)^0=(Z_G(s)\cap G^\vt)^0=Z_{G^{\vt0}}(s)^0$.) From Lemma 4.5 we see also that $W^D$ is
the same as $W_{Z_G(\vt)^0}=W_{Z_{G^{\vt0}}}$. (We use 1.8(a)). From 1.11 we now see that

(a) $\cs_2(G,D)=\cup_{s\in X}j_{W_{Z_{G^{\vt0}}(s)^0}}^{W^D}(\cs_1(Z_G(s),D_s))$
\nl
where $W_{Z_{G^{\vt0}}(s)^0}$ is viewed as a subgroup of $W_{Z_{G^{\vt0}}}=W^D$ as in 1.6(c).
As in 4.7(b), the union in (a) can be restricted to a subset $X'$ of $X$, namely:

(b) $\cs_2(G,D)=\cup_{s\in X'}j_{W_{Z_{G^{\vt0}}(s)^0}}^{W^D}(\cs_1(Z_G(s),D_s))$
\nl
where $X'$ consists of those $s\in X$ such that $Z_{G^{\vt0}}(s)^0$ has the same semisimple rank as $G^{\vt0}$.

\subhead 4.9\endsubhead
In this subsection we assume that we are in type ${}^2E_6$ (see 4.1) so that $G^{\vt0}$ is of type $F_4$.
We write the set $\bI$ in 4.1 as $\{1,2,3,4\}$ in such a way that $\b_1+\b_2$, $\b_2+\b_3$, $\b_2+2\b_3$,
$\b_3+\b_4$ are roots of $G^{\vt0}$ with respect to $T'$; the corresponding simple reflections in
$W^D=W_{G^{\vt0}}$ are denoted by $s_1,s_2,s_3,s_4$ respectively. Let $\o:W^D@>>>W^D$ be the Coxeter group
automorphism given by $s_1\m s_4,s_2\m s_3,s_3\m s_2,s_4\m s_1$. The reflection with respect
to the highest root of $G^{\vt0}$ is the element of $W^D$ given by
$s'_0:=s_1s_2s_3s_2s_1s_4s_3s_2s_3s_4s_1s_2s_3s_2s_1$. Let $s_0=\o(s'_0)\in W^D$.

If $p\ne2$, the one element subsets $K$ in $\cp'(\tI)$ (see 4.1), the corresponding $G_{\tI-K}$ (see 4.4)
and $W_{\tI-K}\sub W^D$ (see 4.7) are as follows:

(i) $K=\{0\}$ with $G_{\tI-K}$ of type $F_4$ with $W_{\tI-K}$ generated by $s_1,s_2,s_3,s_4$;

(ii) $K=\{4\}$ with $G_{\tI-K}$ of type $A_1\T B_3$ with $W_{\tI-K}$ generated by $s_1,s_2,s_3,s_0$;

(iii) $K=\{3\}$ with $G_{\tI-K}$ of type $A_2\T A_2$ with $W_{\tI-K}$ generated by $s_1,s_2,s_4,s_0$ (if
$p\ne3$);

(iv) $K=\{2\}$ with $G_{\tI-K}$ of type $A_3\T A_1$ with $W_{\tI-K}$ generated by $s_1,s_3,s_4,s_0$;

(v) $K=\{1\}$ with $G_{\tI-K}$ of type $C_4$ with $W_{\tI-K}$ generated by $s_2,s_3,s_4,s_0$.
\nl
In each case $\cs_1(G_{\tI-K},G_{\tI-K})\sub\Irr(W_{\tI-K})$ is explicitly known: in case (i), see
\cite{Sh80}; in case (ii),(v) see \cite{L84a}; in case (iii),(iv) it is equal to $\Irr(W_{\tI-K})$.
From this and 4.7(b) we see that: 

(a) {\it if $p\ne2$ then $\cs_2(G,D)$ consists of $1_0,4_1,9_2,12_4,16_5,9_{10},4_{13},1_{24}$; both $8_3,8_3$;
both $8_9,8_9$; one of the two $6_6,6_6$, namely the one of the form $j_{W_{\tI-\{3\}}}^{W^D}(1_6)$;
one of the two $9_6,9_6$, namely the one of the form $j_{W_{\tI-\{2\}}}^{W^D}(1_6)$; one of the two
$4_7,4_7$, namely the one of the form $j_{W_{\tI-\{2\}}}^{W^D}(1_7)$; one of the two $1_{12},1_{12}$,
namely the one of the form $j_{W_{\tI-\{1\}}}^{W^D}(1_{12})$; one of the two $2_{16},2_{16}$, namely the
one of the form $j_{W_{\tI-\{1\}}}^{W^D}(1_{16})$.}
\nl
Here $d_n$ denotes an irreducible representation $E$ of a Weyl group of degree $d$ with $b_E=n$. (For a
Weyl group of type $F_4$ there are at most two irreducible representations with a given $d,n$.)

If $p=2$, the set $X'$ in 4.8 has two elements $s=1$ and $s=g$ where $g\in G^{\vt0}$ has order $3$
with $Z_{G^{\vt0}}(g)^0$ of type $A_2\T A_2$; $Z_G(g)$ is of
type $A_2\T A_2\T A_2$ with $\Ad(\vt)$ interchanging two of these $A_2$-factors;
$W_{Z_{G^{\vt0}}(g)^0}$ can be identified with the subgroup of $W^D$ generated by $s'_0,s_1,s_3,s_4$.

In each case $\cs_1(Z_G(s),D_s)\sub\Irr(W_{Z_{G^{\vt0}}(s)^0})$ is explicitly known: for $s=1$ see
\cite{M05}; for $s=g$ we have $\cs_1(Z_G(s),D_s)=\Irr(W_{Z_{G^{\vt0}}(s)^0})$. From this and 4.8(b) we see
that: 

(b) {\it if $p=2$, then $\cs_2(G,D)$ consists of the same irreducible representations of $W^D$ as those in
(a).}
\nl
Note that

(c) {\it in the union 4.8(b) the term corresponding to $s=g$ is contained in the term corresponding
to $s=1$.}
\nl
Fom (a),(b) we see that

(d) {\it $\cs_2(G,D)$ is independent of $\kk$.}

\subhead 4.10\endsubhead
In this subsection we assume that we are in type ${}^3D_4$ (see 4.1) so that $G^{\vt0}$ is of type $G_2$.
We write the set $\bI$ in 4.1 as $\{1,2\}$ in such a way that $\b_1+\b_2$, $\b_1+2\b_2$, $\b_1+3\b_3$
are roots of $G^{\vt0}$ with respect to $T'$; the corresponding simple reflections in $W^D=W_{G^{\vt0}}$ are
denoted by $s_1,s_2$ respectively. Let $\o:W^D@>>>W^D$ be the Coxeter group
automorphism given by $s_1\m s_2,s_2\m s_1$. The reflection with respect to the highest root of $G^{\vt0}$
is the element of $W^D$ given by $s'_0:=s_1s_2s_1s_2s_1$. Let $s_0=\o(s'_0)\in W^D$.

If $p\ne3$, the one element subsets $K$ in $\cp'(\tI)$ (see 4.1), the corresponding $G_{\tI-K}$ (see 4.4)
and $W_{\tI-K}\sub W^D$ (see 4.7) are as follows:

(i) $K=\{0\}$ with $G_{\tI-K}$ of type $G_2$ with $W_{\tI-K}$ generated by $s_1,s_2$;

(ii) $K=\{2\}$ with $G_{\tI-K}$ of type $A_1\T A_1$ with $W_{\tI-K}$ generated by $s_1,s_0$ (if $p\ne2$)

(iii) $K=\{1\}$ with $G_{\tI-K}$ of type $A_2$ with $W_{\tI-K}$ generated by $s_2,s_0$.
\nl
In each case, $\cs_1(G_{\tI-K},G_{\tI-K})\sub\Irr(W_{\tI-K})$ is explicitly known: in case (i), see
\cite{Sp85}; in case (ii),(iii) it is equal to $\Irr(W_{\tI-K})$. From this and 4.7(b) we see that: 

(a) {\it if $p\ne3$ then $\cs_2(G,D)$ consists of $1_0,2_1,2_2,1_6$;
one of the two $1_3,1_3$, namely the one of the form $j_{W_{\tI-\{1\}}}^{W^D}(1_3)$.}
\nl
(Notation as in 4.9(a).)

If $p=3$, the set $X'$ in 4.8 has two elements $s=1$ and $s=g$ where $g\in G^{\vt0}$ has order $2$
with $Z_{G^{\vt0}}(g)^0$ of type $A_1\T A_1$; $Z_G(g)$ is of type $A_1\T A_1\T A_1\t A_1$ with
$\Ad(\vt)$ permuting cyclically three of these $A_1$-factors;
$W_{Z_{G^{\vt0}}(g)^0}$ can be identified with the subgroup of $W^D$ generated by $s'_0,s_2$.

In each case $\cs_1(Z_G(s),D_s)\sub\Irr(W_{Z_{G^{\vt0}}(s)^0})$ is explicitly known: for $s=1$ see
\cite{M05}; for $s=g$ we have $\cs_1(Z_G(s),D_s)=\Irr(W_{Z_{G^{\vt0}}(s)^0})$. From this and
4.8(b) we see that:

(b) {\it if $p=3$, then $\cs_2(G,D)$ consists of the same irreducible representations of $W^D$
as those in (a).}
\nl
Note that

(c) {\it in the union 4.8(b), the term corresponding to $s=g$
is contained in the term corresponding to $s=1$.}
\nl
Fom (a),(b) we see that

(d) {\it $\cs_2(G,D)$ is independent of $\kk$.}

\subhead 4.11\endsubhead
In this subsection we assume that we are in type ${}^2A_{2n+1},n\ge1$, or type ${}^2A_{2n},n\ge1$, or type
${}^2D_n,n\ge4$. If $p\ne2$, the Strata of $D$ and the set $\cs_2(G,D)$ are determined by 3.1(a), 3.3(a) (we
use also 1.13); we see that $\cs_2(G,D)$ can be identified with $BP^n_{3,1},BP^n_{1,3},BP^{n-1}_{4,0}$
respectively. If $p=2$, from 4.8(b) we have

(a) $\cs_2(G,D)=\cs_1(G,D)$
\nl
and this is identified in \cite{L04,\S13},
\cite{MS04} with ${}^m\hBP^n_{3,1}$, ${}^m\hBP'{}^n_{1,3}$, ${}^m\hBP^{n-1}_{4,0}$ respectively (with $m\gg0$).
(The proof of this identification in {\it loc.cit.} is very similar to the analogous statement for connected
classical groups in characteristic $2$ given in \cite{LS85}.) Using the identifications in 2.2, we see that
$\cs_2(G,D)$ is identified with $BP^n_{3,1}$, $BP^n_{1,3}$, $BP^{n-1}_{4,0}$ respectively. We now see that

(b) {\it $\cs_2(G,D)$ is independent of $\kk$.}

\subhead 4.12\endsubhead
We prove Theorem 1.15. By a sequence of reductions as in \cite{L04, 12.1-12.7} we see that we can assume that
$G,D,\d$ are as in the first sentence of 4.1. If $\d=1$, we have $G=G^0=D$ and the desired result follows
from \cite{L15}. Assume now that $\d\ge2$. In this case the desired result follows from 4.9(d), 4.10(d),
4.11(b). The theorem is proved.

\head 5. Further results\endhead
\subhead 5.1\endsubhead
In this subsection we assume that $D=G^0=G$ and $\kk=\CC$.
Let $\sgn_W$ be the sign representation of $W$. The following is well known:

(a) {\it If $g\in G$ is unipotent then $E_g=\sgn_W$ if and only if $g=1$.}
\nl
Let ${}^0\cs_2(W)$ be the set of all $E\in\Irr(W)$ such that $E=j_{W'}^W(\sgn_{W'})$ where $W'$ is the Weyl
group of $Z_G(s)^0$ (for some semisimple $s\in G$) viewed as a subgroup of $W$.
We have ${}^0\cs_2(W)\sub\cs_2(W)$ (we use (a)). Note that ${}^0\cs_2(W)$ parametrizes
a subset of the set of unipotent classes of the group of type dual to that of $G$.

Let $\Si\in Str(G)$ be corresponding to $E\in\cs_2(W)$. From (a) and the definitions we deduce:

(b) {\it $\Si$ contains some semisimple element if and only if $E\in{}^0\cs_2(W)$.}

\subhead 5.2\endsubhead
We preserve the setup of 5.1.
Let $G_0$ be a maximal compact subgroup of $G$. Then $G_0$ is partitioned into subsets defined by the
intersections of the various Strata of $G$ with $G_0$. The Strata of $G$ which have nonempty intersection
with $G_0$ are precisely those in 5.1(b).

\subhead 5.3\endsubhead
We preserve the setup of 5.1. For $a\in\NN$ we set

$Y_a=(0,1,0,1,0,1,\do,0,1,0,0,0,\do)\in BP^a_{1,0}$,

$Y'_a=(1,0,1,0,1,\do,0,1,0,0,0,\do)\in BP^a_{0,1}$.
\nl
Assume that $G=Sp_{2n}(\kk)$, $n\ge2$. We identify $\cs_2(W)$ with $BP^n_{2,2}$ as in \cite{L15, 3.5(b)}.
Then ${}^0\cs_2(W)$ becomes the set of bipartitions of the form $Y_a+Y'_b+C$ where $C\in BP^{n-a-b}_{0,0}$,
$a+b\le n$.

Assume that $G=SO_{2n+1}(\kk)$, $n\ge2$. We identify $\cs_2(W)$ with $BP^n_{2,2}$ as in \cite{L15, 3.5(b)}.
Then ${}^0\cs_2(W)$ becomes the set of bipartitions of the form $Y_a+Y_b+C$ where $C\in BP^{n-a-b}_{0,0}$,
$a+b\le n$.

Assume that $G=SO_{2n}(\kk)$, $n\ge4$. We identify $\cs_2(W)$ (modulo the action of $O_{2n}(\kk)$)
with $BP^n_{0,4}$ as in \cite{L15, 3.10(b)}. Then ${}^0\cs_2(W)$ (modulo the action of $O_{2n}(\kk)$)
becomes the set of bipartitions of the form $Y'_a+Y'_b+C$ where $C\in BP^{n-a-b}_{0,0}$, $a+b\le n$.

\subhead 5.4\endsubhead
We no longer assume that $D=G^0=G$ but we assume that $p=0$.
Let $s_0\in D$ be a semisimple element and let $T_0$ be a maximal torus of $Z_G(s_0)^0$.
It is known that $T_0s_0$ consists of semisimple elements and any semisimple element in $D$ is
$G^0$-conjugate to an element of $T_0s_0$ (see \cite{L03, 1.14}). In this case the description of $\cs_2(G,D)$
in 1.11 simplifies as follows: the representations of $W^D$ in $\cs_2(G,D)$ are exactly those
obtained by applying
$j_{W_s}^{W^D}$ to any representation of $W_s$ in $\cs_1(Z_G(s)^0,Z_G(s)^0)$ where $s$ is any 
element in $T_0s_0$; here $W_s=W_{Z_G(s)^0}$ is viewed as a subgroup of $W^D$ as in 1.6(c).
(We use that $D_1$ in 1.11 is now $Z_G(s)^0$.) For $s,W_s$ above we consider a point $\x$ on the torus
$T_0^*$ dual to $T_0$ and we denote by $W_{s,\x}$ the subgroup of $W_s$
generated by the reflections which keep $\x$ fixed (in the natural action of $W_s$ on $T_0^*$). 
From \cite{L09} it is known that the representations of $W_s$ in $\cs_1(Z_G(s)^0,Z_G(s)^0)$ are exactly those
obtained by applying $j_{W_{s,\x}}^{W_s}$ to any special representation of $W_{s,\x}$. Let $\cp(G,D)$ be the
collection of reflection subgroups of $W^D$ which are conjugate to a subgroup  of the form $W_{s,\x}$ for some
$s\in T_0s_0$ and some $\x\in T_0^*$; the subgroups in $\cp(G,D)$ are said to be the $2$-parabolic
subgroups of $W^D$. It follows that the representations of $W^D$ in $\cs_2(G,D)$ are exactly
those obtained by applying $j_{W_1}^{W^D}$ (with $W_1\in\cp(G,D)$) to any special representation of $W_1$.

The set $\cp(G,D)$ consists of the reflection subgroups of $W^D$ of type

(i) $B_a\T B_b\T B_c\T B_d\T\ca$    if $(G,D)$ is of type ${}^2D_{n+1},n\ge3$ (see 4.1);

(ii) $B_a\T B_b\T B_c\T D_d\T\ca$   if $(G,D)$ is of type ${}^2A_{2n},n\ge1$, (see 4.1);

(iii) $B_a\T B_b\T D_c\T D_d\T\ca$   if $G=G^0=D$ is $Sp_{2n}(\kk)$ or $SO_{2n+1}(\kk)$, $n\ge2$;

(iv) $B_a\T D_b\T D_c\T D_d\T\ca$   if $(G,D)$ is of type ${}^2A_{2n-1},n\ge2$, (see 4.1);

(v) $D_a\T D_b\T D_c\T D_d\T\ca$   if $G=G^0=D$ is $SO_{2n}(\kk),n\ge4$.
\nl
Note that in (i)-(iv), $W^D$ is of type $B_n$ while in (v), $W^D=W$ is of type $D_n$; in each case we have
$(a,b,c,d)\in\NN^4$, $a+b+c+d\le n$; $\ca$ stands for a product of reflection subgroups of type $A$. 

\subhead 5.5\endsubhead
In this subsection we assume that $G,D,\d$ are as in 4.1 with $p=\d\in\{2,3\}$. We show:

(a) {\it Each Stratum of $D$ contains a unique unipotent $G^0$-orbit.}
\nl
It is enough to show that $\cs_2(G,D)=\cs_1(G,D)$. If $G,D,\d$ is as in 4.11 this follows from 4.11(a). If
$G,D,\d$ is as in 4.9 (resp. 4.10) this follows from 4.9(c) (resp. 4.10(c)).

\subhead 5.6\endsubhead
Let $\cs(W)$ be the subset of $\Irr(W)$ consisting of special representations.
We define a map $\z:Str(D)@>>>\cs(W)$ as the composition
$Str(D)=\cs_2(W,[D])\sub\Irr(W)@>\ti\z>>\cs(W)$ where $\ti\z$ is defined as follows.
If $D=G^0$, $\ti\z:\Irr(W)@>>>\cs(W)$ associates to any $E\in\Irr(W)$ the special representation in the same
family as $E$. In the general case, we can find a connected
reductive group $\cg$ defined over $F_q$ with connected centre, with Weyl group
$W$ and with Frobenius map acting on $W$ as $[D]$. Then $\Irr(W^D)$ is in natural bijection with the
set of irreducible representations of $\cg(F_q)$ which have nonzero invariants under a Borel subgroup
defined over $F_q$; this is a subset of the set of unipotent representations of $\cg(F_q)$ which by its
known parametrization \cite{L84b}
is decomposed into pieces indexed by the families of $W$ which are stable under $[D]$;
in this way we can associate to each representation in $\Irr(W^D)$ a family of $W$ and hence the unique
special representation of $W$ in that family. This defines $\ti\z$ and hence $\z$.

\subhead 5.7\endsubhead
We define a partition of $D$ into ``special pieces'' indexed by $\cs(W)$. For $E\in\cs(W)$, the special piece
corresponding to $E$ is $\cup_{\Si\in Str(D);\z(\Si)=E}\Si$. We expect
that this special piece is locally closed in $D$. As supporting evidence, we note that, when $G=G^0=D,p=0$, the
intersection of this special piece with the unipotent variety of $G$ coincides with a special piece of that
unipotent variety considered in \cite{L97}; in particular it is locally closed in the unipotent variety.

\subhead 5.8\endsubhead
In this subsection we assume that that $\kk$ is an algebraic closure of a finite field $F_q$
and that $G=G^0=D$ has connected centre; we also assume that $G$ has a fixed $F_q$-rational structure.
Let $G^*$ be a connected reductive group over
$\kk$ of type dual to that of $G$ with the induced $F_q$-structure.
Let $sc(G^*)$ be the set of all special conjugacy classes in $G^*$ which are defined over $F_q$. It is known
\cite{L84b, 13.2} that there is a natural map $\t:\Irr(G(F_q))@>>>sc(G^*)$. For any $\Si\in Str(G^*)$ let
$\Irr(G(F_q))_\Si$ be the set of all $\EE\in\Irr(G(F_q))$ such that $\t(\EE)\sub\Si$. Thus we have a partition

(a) $\Irr(G(F_q))=\sqc_{\Si\in Str(G^*)}\Irr(G(F_q))_\Si$.
\nl
We have natural bijections  $a:Str(G^*)@>\si>>Str(G)$, $b:Str(G)@>\si>>Str(G_\CC)$
where $G_\CC$ denotes a connected reductive group over $\CC$ of the same type as that of $G$ (see \cite{L15}:
each of the sets $Str(G^*),Str(G),Str(G_\CC)$ is indexed by $\cs_2(W)$.
From the results in \cite{L84b, 13.3}, we see that for $\Si\in Str(G^*)$, $\Irr(G(F_q))_\Si$ is empty unless
$ba(\Si)$ contains some unipotent class of $G_\CC$.
Thus (a) can be regarded as a partition of $\Irr(G(F_q))$ into pieces indexed by the Strata $\Si'$
of $G$ such that $b(\Si')$ contains a unipotent class in
$G_\CC$. (In {\it loc.cit.} this partition is defined only under the assumption that $p$ is a good prime for
$G$ but the same definition applies without this assumption; if $p$ is a bad prime for $G$, not all $\Si'$ as
above contribute to the partition.)

\enddocument